\documentclass{pnastwo}

\usepackage{amssymb,amsfonts,amsmath}

\url{www.pnas.org/cgi/doi/10.1073/pnas.1313071111}
\copyrightyear{2014}
\issuedate{July 8, 2014}
\volume{vol. 111}
\issuenumber{no. 27}

\newtheorem{example} {Example}

\def \Rset {{\mathbb R}}         
\def \Cset {{\mathbb C}}
\def \KK {{\mathbb K}}

\def \Zset {{\mathbb Z}}

\def \TT {{\mathcal{T}}}

\def \UU {{\mathcal{U}}}

\def \TT {{\mathcal{T}}} 

\def \qb {{\bf{q}}}

\def \ex {{\bf{ex}}}

\def \de {\delta}
\def \al {\alpha}
\def \be {\beta}
\def \vpi {\varpi}
\def \la {\lambda}

\def \Om {\Omega}
\def \ga {\gamma}
\def \de {\delta}

\def \sig {\sigma}

\def \vp {\varphi}
\def \ep {\epsilon}
\def \De {\Delta}
\def \sig{\sigma}

\def \mt  {\mapsto}

\def \bu  {\bullet}                 

\def \rcor {\rangle}
\def \lcor {\langle}

\def \ol {\overline}
\def \wt {\widetilde}

\def \id { {\mathrm{id}} }

\def \sign { {\mathrm{sign}} }

\DeclareMathOperator \rk { {\mathrm{rk}} }
\DeclareMathOperator \range { {\mathrm{range}} }

\def \g  {\mathfrak{g}}   

\def \sl {\mathfrak{sl}}

\def \n  {\mathfrak{n}}

\def \sl {\mathfrak{sl}}

\DeclareMathOperator \Span { {\mathrm{Span}} }

\DeclareMathOperator \Fract { {\mathrm{Fract}} }

\renewcommand \max { {\mathrm{max}} }
\newcommand\kx{\KK^*}

\newcommand\HH{{\mathcal{H}}}

\begin{document}



\title{Quantum cluster algebras and
quantum nilpotent algebras}





\author{K. R. Goodearl\affil{1}{Department of Mathematics, 
University of California, Santa Barbara, CA 93106, U.S.A.}
\and M. T. Yakimov\affil{2}{Department of Mathematics, 
Louisiana State University, Baton Rouge, LA 70803, U.S.A.}
}

\contributor{Proceedings of the National Academy of Sciences
of the United States of America 111 (2014) 9696--9703}

\maketitle

\begin{article}

\begin{abstract} 
A major direction in the theory of cluster algebras is to construct
(quantum) cluster algebra structures on the (quantized) coordinate 
rings of various families of varieties arising in Lie theory. 
We prove that all algebras in a very large axiomatically defined class of 
noncommutative algebras possess canonical quantum cluster algebra structures.
Furthermore, they coincide with the corresponding 
upper quantum cluster algebras. We also establish analogs of these results 
for a large class of Poisson nilpotent algebras. Many important families of 
coordinate rings are subsumed in the class we are covering, which leads to a 
broad range of application of the general results to the above mentioned 
types of problems. As a consequence, we prove the Berenstein--Zelevinsky conjecture for the 
quantized coordinate rings of double Bruhat cells and construct 
quantum cluster algebra structures on all quantum unipotent groups, 
extending the theorem of Gei\ss, Leclerc and Schr\"oer for the case of 
symmetric Kac--Moody groups. Moreover, we prove that the upper 
cluster algebras of Berenstein, Fomin and Zelevinsky associated 
to double Bruhat cells coincide with the corresponding cluster 
algebras. 
\end{abstract}






\dropcap{T}he theory of cluster algebras provides a unified framework 
for treating a number of problems in diverse areas of mathematics such as 
combinatorics, representation theory, topology, mathematical physics,
algebraic and Poisson geometry, and dynamical systems \cite{FZ2,FG,KS,GSV-b,K,W,GLSr}. 
The construction of cluster algebras was invented 
by Fomin and Zelevinsky \cite{FZ2} who also obtained a number 
of fundamental results on them. This construction 
builds algebras in a novel way by producing infinite generating 
sets via a process of mutation rather than the classical 
approach using generators and relations.

The main algebraic approach to cluster algebras 
relies on representations of finite dimensional 
algebras and derived categories \cite{K,R}. In this paper 
we describe a different algebraic approach based on noncommutative 
ring theory.

An important range of problems in the theory of cluster algebras is
to prove that the coordinate rings of certain algebraic varieties 
coming from Lie theory admit cluster algebra structures.
The idea is that, once this is done, one can use cluster algebras 
to study canonical bases in such coordinate rings.
Analogous problems deal with the corresponding quantizations.
The approach via representation theory to this type of problem
needs combinatorial data for quivers as a starting point. 
Such might not be available a priori. This approach also differs  
from one family of varieties to another, which means that in each case 
one needs to design an appropriate categorification process.

We prove that all algebras in a very general, 
axiomatically defined class of quantum nilpotent algebras 
are quantum cluster algebras. 
The proof is based on 
constructing quantum clusters by 
considering sequences of prime elements in chains of 
subalgebras which are noncommutative unique factorization 
domains. 
These clusters are canonical relative to the mentioned chains of subalgebras, which are determined by the presentation of the quantum nilpotent algebra.
The construction does not rely on any initial 
combinatorics of the algebras. On the contrary, the construction itself
produces intricate combinatorial data for prime elements in chains 
of subalgebras. When this is applied to special cases, we recover 
the Weyl group combinatorics which played a key role 
in categorification earlier \cite{GLSr,FZ1,BFZ}. Because of this, we expect that our construction will be 
helpful in building a unified categorification of quantum nilpotent algebras.
Finally, we also prove similar results for (commutative) cluster algebras 
using Poisson prime elements. 

The results of the paper have many applications since 
a number of important families of algebras arise as special 
cases in our axiomatics. Berenstein, Fomin, and Zelevinsky 
proved \cite{BFZ} that the coordinate rings of all 
double Bruhat cells in every complex simple Lie group are 
upper cluster algebras. It was an important problem 
to decide whether the latter coincide with the corresponding 
cluster algebras. We resolve this problem positively. 
On the quantum side, we prove the Berenstein--Zelevinsky
conjecture \cite{BZ} on cluster algebra structures for all 
quantum double Bruhat cells. Finally, we establish
that the quantum Schubert cell algebras for all complex simple Lie groups 
have quantum cluster algebra structures. Previously 
this was known for symmetric Kac--Moody groups 
due to Gei\ss, Leclerc and Schr\"oer \cite{GLS}. 
\section{Prime elements of quantum nilpotent algebras}

\subsection{Definition of quantum nilpotent algebras}
Let $\KK$ be an arbitrary base field.
A skew polynomial extension of a $\KK$-algebra $A$,
$$
A \mt A[x; \sig, \de],
$$
is a generalization of the classical polynomial algebra $A[x]$.
It is defined using an algebra automorphism $\sig$ of $A$ and 
a skew-derivation $\de$. The algebra $A[x; \sig, \de]$ is isomorphic to 
$A[x]$ as a vector space and the variable $x$ commutes with 
the elements of $A$ as follows:
$$
x a = \sig(a) x + \de(a) \quad {\mbox{for all}} \; \;  a \in A.
$$

For every nilpotent Lie algebra $\n$ of dimension $m$
there exists a chain of ideals of $\n$
$$
\n = \n_m  \vartriangleright \n_{m-1} \vartriangleright \ldots \vartriangleright \n_1 \vartriangleright
\n_0 = \{ 0 \}
$$
such that $\dim (\n_k /\n_{k-1}) =1$ and $[\n,\n_k] \subseteq \n_{k-1}$, 
for $1 \leq k \leq m$.
Choosing an element $x_k$ in the complement of $\n_{k-1}$ in $\n_k$ 
for each $1 \leq k \leq m$
leads to 
the following presentation of the universal enveloping 
algebra $\UU(\n)$ as an iterated skew polynomial extension: 
$$
\UU(\n) \cong \KK [x_1][x_2; \id, \de_2] \ldots [x_m; \id, \de_m]
$$
where all the derivations $\de_2, \ldots, \de_m$ are locally nilpotent.

\begin{definition}
\label{q-nil}
An iterated skew polynomial extension 
\begin{equation}
\label{R-quant-nilp}
R = \KK[x_1][x_2; \sig_2, \delta_2] \cdots [x_N; \sig_N, \delta_N]
\end{equation}
is called a quantum nilpotent algebra if it is 
equipped with a rational action of a $\KK$-torus $\HH$ 
by $\KK$-algebra automorphisms satisfying the following conditions:

{\rm(a)} The elements $x_1,\dots,x_N$ are $\HH$-eigenvectors.

{\rm(b)} For every $2 \leq k \leq N$, $\de_k$ is a locally nilpotent 
$\sig_k$-der\-i\-va\-tion of $\KK[x_1] \cdots [x_{k-1}; \sig_{k-1}, \delta_{k-1}]$. 

{\rm(c)}  For every $1 \leq k \leq N$, there exists $h_k \in \HH$ such that 
$\sig_k = (h_k \cdot)$ and the $h_k$-eigenvalue of $x_k$, to be denoted by $\la_k$, is not a root 
of unity.
\end{definition}

The universal enveloping algebras of finite dimensional 
nilpotent Lie algebras satisfy all of the conditions 
in the definition except for the last part of the third one. More precisely, 
in that case one 
can take $\HH = \{1\}$, conditions (a)--(b) are satisfied,
and in condition (c) we have $\la_k=1$. Thus,
condition (c) is the main feature that separates the class 
of quantum nilpotent algebras from the class of universal 
enveloping algebras of nilpotent Lie algebras. The torus 
$\HH$ is needed in order to define the eigenvalues 
$\la_k$. 

The algebras in Definition \ref{q-nil} are also known as 
Cauchon--Goodearl--Letzter (CGL) extensions. The axiomatics came from
the works \cite{GL,Ca1} which investigated in this generality 
the stratification of the prime spectrum of an algebra into strata 
associated to its $\HH$-prime ideals.

The Gelfand--Kirillov dimension of the algebra in Eq.~\eqref{R-quant-nilp} equals $N$. 

\begin{example} 
\label{quant-matr}
For two positive integers $m$ and $n$, and
$q \in \kx$, define the algebra of quantum matrices $R_q[M_{m \times n}]$ 
as the algebra with generators $t_{ij}$, $1 \leq i \leq m$ and 
$1 \leq j \leq n$, and relations 
\begin{align*}
t_{ij} t_{kj} &= q t_{kj} t_{ij}, &&\mbox{for} \; \; i < k, \\
t_{ij} t_{il} &= q t_{il} t_{ij}, &&\mbox{for} \; \; j < l, \\ 
t_{ij} t_{kl} &= t_{kl} t_{ij}, &&\mbox{for} \; \; i < k, \; j > l,\\
t_{ij} t_{kl} - t_{kl} t_{ij} &= (q-q^{-1}) t_{il} t_{kj}, &&\mbox{for} \; \; i < k, \; j<l.
\end{align*}
It is an iterated skew polynomial extension where
$$
R_q[M_{m \times n}] = \KK[x_1][x_2; \sig_2, \de_2] \ldots [x_N; \sig_N, \de_N], 
$$
$N=mn$, and $x_{(i-1)n +j} = t_{ij}$. It is easy to write explicit formulas 
for the automorphisms $\sig_k$ and skew derivations $\de_k$ from 
the above commutation relations, and to check that each $\de_k$ is locally nilpotent.
The torus $\HH = (\kx)^{m+n}$ acts on $R_q[M_{m \times n}]$
by algebra automorphisms by the rule
$$
(\xi_1, \ldots, \xi_{m+n}) \cdot t_{ij} := \xi_i \xi_{m+j}^{-1} t_{ij}
$$
for all $(\xi_1, \ldots, \xi_{m+n}) \in (\kx)^{m+n}$. Define 
$$
h_{ij} := (1, \ldots, 1, q^{-1}, 1, \ldots, 1, q,  1, \ldots, 1) \in \HH
$$
where $q^{-1}$ and $q$ reside in positions $i$ and $m+j$, respectively.
Then $\sig_{(i-1)n+j} = (h_{ij} \cdot)$ and 
$$
h_{ij} \cdot t_{ij} = q^{-2} t_{ij}.
$$
Thus, for all $q \in \kx$ which are not roots of unity, the algebras 
$R_q[M_{m \times n}]$ are examples of quantum nilpotent algebras.
\end{example}

\subsection{Unique factorization domains}
The notion of unique factorization domain plays an important role in 
algebra and number theory. Its noncommutative analog was introduced
by Chatters in \cite{Cha}. A nonzero, non-invertible element $p$ of a domain $R$ (a ring without zero 
divisors) 
is called prime if $p R = Rp$ and the factor $R/Rp$ is a domain.
A noetherian domain $R$ is called a unique factorization 
domain (UFD) if every nonzero prime ideal of $R$ contains a 
prime element. Such rings possess the unique factorization 
property for all of their nonzero normal elements -- the 
elements $u \in R$ with the property that $Ru = uR$. If the ring $R$ is acted 
upon by a group $G$, then one can introduce an equivariant 
version of this property: Such an $R$ is called a $G$-UFD
if every nonzero $G$-invariant prime ideal of $R$ contains a prime 
element which is a $G$-eigenvector.

It was shown in \cite{LLR} that every quantum nilpotent algebra
$R$ is a noetherian $\HH$-UFD. An $\HH$-eigenvector of such a 
ring $R$ will be called a homogeneous element since 
this corresponds to the homogeneity property 
with respect to the canonical induced grading of $R$ 
by the character lattice of $\HH$. In particular, we will 
use the more compact term of homogeneous prime element 
of $R$ instead of a prime element of $R$
which is an $\HH$-eigenvector.

\subsection{Sequences of prime elements} 
Next, we classify the 
set of all homogeneous prime elements 
of the chain of subalgebras
$$
\{ 0 \} \subset 
R_1 \subset R_2 \subset \ldots \subset R_N = R
$$
of a quantum nilpotent algebra $R$,
where $R_k$ is the subalgebra of $R$ generated by 
the first $k$ variables $x_1, \ldots, x_k$.

We will denote by $\Zset$ and $\Zset_{\geq 0}$ the sets of all integers 
and nonnegative integers, respectively.
Given two integers $l \leq  k$, set $[l,k] :=\{ l, l+1, \ldots, k\}$. 

For a function $\eta : [1,N] \to \Zset$, define the predecessor function
$p : [1,N] \to [1,N] \sqcup \{ - \infty \}$ and 
successor function
$s : [1,N] \to [1,N] \sqcup \{ + \infty \}$ for its level sets
by
$$
p(k) = 
\begin{cases}
\max \{ l <k \mid \eta(l) = \eta(k) \}, 
&\text{if such $l$ exists}, 
\\
- \infty, \; & \text{otherwise}, 
\end{cases}
$$
and 
$$
s(k) = 
\begin{cases}
\min \{ l > k \mid \eta(l) = \eta(k) \}, 
&\text{if such $l$ exists}, 
\\
+ \infty, \; & \text{otherwise}. 
\end{cases}
$$

\begin{theorem}
\label{chain-prime}
Let $R$ be a quantum nilpotent algebra of dimension $N$.
There exist a function $\eta : [1,N] \to \Zset$ 
and elements
$$
c_k \in R_{k-1} \; \;  \mbox{for all} \; \; 2 \leq k \leq N \; \; 
\mbox{with} \; \; \de_k \neq 0
$$
such that the elements $y_1, \ldots, y_N \in R$, recursively defined by 
$$
y_k := 
\begin{cases}
y_{p(k)} x_k - c_k, &\text{if} \; \;  \de_k \neq 0 \\
x_k, & \text{if} \; \; \de_k = 0,  
\end{cases}
$$
are homogeneous and have the property that for every $k \in [1,N]$
the homogeneous prime elements of $R_k$ are precisely the nonzero 
scalar multiples of the elements
$$
y_l \; \; \mbox{for} \; \;  l \in [1,k] \; \; \text{with} \; \; s(l) > k .
$$ 
In particular, $y_k$ is a homogeneous prime element of $R_k$, for all $k \in [1,N]$.
The sequence $y_1, \ldots, y_N$ and the level sets 
of a function $\eta$ with these properties are both unique.
\end{theorem}

\begin{example}
\label{y-quant-matr} Given two subsets $I = \{ i_1 < \cdots < i_k \} \subset [1,m]$
and $J = \{ j_1 < \cdots < j_k \} \subset [1,n]$,
define the quantum minor $\De_{I,J} \in R_q[M_{m \times n}]$
by 
$$
\De_{I,J} = \sum_{\sig \in S_k} (-q)^{\ell(\sig)}
t_{i_1 j_{\sig(1)}} \ldots t_{i_k j_{\sig(k)}}
$$
where $S_k$ denotes the symmetric group and $\ell : S_k \to \Zset_{\geq 0}$ 
the standard length function.

For the algebra of quantum matrices $R_q[M_{m \times n}]$, 
the sequence of prime elements from Theorem {\rm\ref{chain-prime}}
consists of solid quantum minors; more precisely,
$$
y_{(i-1)n+j} = \De_{[i- \min(i,j)+1, i], [j- \min(i,j) + 1, j]}
$$
for all $1 \leq i \leq m$ and $1 \leq j \leq n$.
Furthermore, the function $\eta : [1,mn] \to \Zset$ can be chosen
as 
$$
\eta((i-1)n +j) := j-i.
$$
\end{example}

\begin{definition}
\label{range} The cardinality of the range of the function $\eta$ 
from Theorem {\rm\ref{chain-prime}} is called the rank of the quantum nilpotent 
algebra $R$ and is denoted by $\rk(R)$.
\end{definition}
For example, the algebra of quantum matrices $R_q[M_{m \times n}]$ has rank 
$m+n -1$.

\subsection{Embedded quantum tori}
An $N \times N$ matrix $\qb:=(q_{kl})$ with entries in $\KK$ is called 
multiplicatively skewsymmetric if 
$$
q_{kl} q_{lk} = q_{kk} = 1 \quad {\mbox{for}} \; \;  1 \leq l,k \leq N.
$$
Such a matrix gives rise to the quantum torus $\TT_\qb$ 
which is the $\KK$-algebra with generators $Y_1^{ \pm 1}, \ldots, Y_N^{ \pm 1}$
and relations 
$$
Y_k Y_l = q_{kl} Y_l Y_k \quad {\mbox{for}} \; \;  1 \leq l,k \leq N.
$$

Let $\{e_1, \ldots, e_N\}$ be the standard basis of the lattice $\Zset^N$.
For a quantum nilpotent algebra $R$ of dimension $N$, 
define the eigenvalues $\la_{kl} \in \KK$:
\begin{equation}
\label{eig}
h_k \cdot x_l = \la_{kl} x_l \quad {\mbox{for}} \; \; 1 \leq l < k \leq N.
\end{equation}
There exists 
a unique group bicharacter $\Om : \Zset^N \times \Zset^N \to \kx$ such that
$$
\Om(e_k, e_l)=
\begin{cases}
1, &\text{if $\; k=l$},
\\
\la_{kl},
&\text{if $\; k>l$}, 
\\
\la_{lk}^{-1},
&\text{if $\; k<l$}.
\end{cases}
$$
Set $e_{- \infty} :=0$.
Define the vectors 
\begin{equation}
\label{e-ol-vect}
\ol{e}_k := e_k + e_{p(k)} + \cdots \in \Zset^N,
\end{equation}
noting that only finitely many terms in the sum are 
nonzero.
Then $\{ \ol{e}_1, \ldots, \ol{e}_N \}$ is another 
basis of $\Zset^N$.

\begin{theorem}
\label{qu-torus-embedding} For each quantum nilpotent 
algebra $R$, the sequence of prime elements from 
Theorem {\rm\ref{chain-prime}} defines an embedding of the 
quantum torus $\TT_\qb$ associated to the $N \times N$ 
multiplicatively skewsymmetric matrix with entries
$$
q_{kl} := \Om(\ol{e}_k, \ol{e}_l), \quad 1 \leq k,l \leq N
$$
into the division ring of fractions $\Fract(R)$ of $R$
such that $Y_k^{\pm 1} \mt y_k^{\pm 1}$, for all $1 \le k \le N$. 
\end{theorem}

\section{Cluster structures on quantum nilpotent algebras}

\subsection{Symmetric quantum nilpotent algebras}
\begin{definition}
\label{symmetric} A quantum nilpotent algebra $R$ as in Definition 
{\rm\ref{q-nil}} will be called symmetric if it can be presented 
as an iterated skew polynomial extension for the reverse order of
its generators,
$$
R = \KK[x_N][x_{N-1}; \sig^*_{N-1}, \delta^*_{N-1}] \cdots [x_1; \sig^*_1, \delta^*_1],
$$
in such a way that conditions {\rm(a)--(c)} in Definition {\rm\ref{q-nil}}
are satisfied for some choice of $h_N^*, \ldots, h_{1}^* \in \HH$.
\end{definition}

A quantum nilpotent algebra $R$ is symmetric if and only if it satisfies 
the Levendorskii--Soibelman type straightening law
$$
x_k x_l - \la_{kl} x_l x_k = 
\sum_{n_{l+1}, \ldots, n_{k-1} \in \Zset_{\geq 0}} 
\xi_{n_{l+1}, \ldots, n_{k-1}} x_{l+1}^{n_{l+1}} \ldots x_{k-1}^{n_{k-1}}
$$
for all $l <k$ (where the $\xi_\bu$ are scalars)   
and there exist $h_k^* \in \HH$ such that $h_k^* \cdot x_l = \la_{lk}^{-1} x_l$
for all $l >k$. The defining commutation relations for the algebras of quantum matrices 
$R_q[M_{m \times n}]$ imply that they are symmetric quantum nilpotent algebras.

Denote by $\Xi_N$ the subset of the symmetric group $S_N$
consisting of all permutations that have the property that
\begin{align}
\label{1}
&\tau(k) = \max \, \tau( [1,k-1]) +1 \quad \mbox{or} 
\\
\label{2}
&\tau(k) = \min \, \tau( [1,k-1]) - 1
\end{align}
for all $2 \leq k \leq N$. In other words, $\Xi_N$ consists of those 
$\tau \in S_N$ such that $\tau([1,k])$ is an interval for all $2 \leq k \leq N$. 
For each $\tau \in \Xi_N$, a symmetric quantum nilpotent algebra $R$ has the 
presentation
\begin{equation}
\label{tau}
R = \KK [x_{\tau(1)}] [x_{\tau(2)}; \sig''_{\tau(2)}, \de''_{\tau(2)}] 
\cdots [x_{\tau(N)}; \sig''_{\tau(N)}, \de''_{\tau(N)}]
\end{equation}
where $\sig''_{\tau(k)} := \sig_{\tau(k)}$ and 
$\de''_{\tau(k)} := \de_{\tau(k)}$ if Eq.~\eqref{1} is satisfied, while 
$\sig''_{\tau(k)} := \sig^*_{\tau(k)}$ and 
$\de''_{\tau(k)} := \de^*_{\tau(k)}$ if 
Eq.~\eqref{2} holds. This presentation satisfies 
the conditions (a)--(c) in Definition \ref{q-nil} for the 
elements $h''_{\tau(k)} \in \HH$, 
given by $h''_{\tau(k)} := h_{\tau(k)}$ in case of Eq.~\eqref{1} 
and $h''_{\tau(k)} := h^*_{\tau(k)}$ in case of Eq.~\eqref{2}.
The use of the term symmetric in Definition \ref{symmetric} is motivated 
by the presentations in Eq.~\eqref{tau} parametrized by the elements
of the subset $\Xi_N$ of the symmetric group $S_N$.

\begin{proposition} 
\label{la-isi}
For every symmetric quantum nilpotent algebra $R$, 
the $h^*_k$-eigenvalues of $x_k$, to be denoted by $\la^*_k$,
satisfy
$$
\la_k^* = \la_l^*
$$
for all $1 \leq k,l \leq N$ such that $\eta(k) = \eta(l)$
and $s(k) \neq + \infty$, $s(l) \neq + \infty$.
They are related to the eigenvalues $\la_l$ by 
$$
\la_k^* = \la_l^{-1}
$$
for all $1 \leq k,l \leq N$ such that $\eta(k) = \eta(l)$
and $s(k) \neq + \infty$, $p(l) \neq - \infty$.
\end{proposition}

\subsection{Construction of exchange matrices}
\  Our construction of a quantum cluster algebra structure 
on a symmetric quantum nilpotent algebra $R$ of dimension $N$ 
will have as the set of exchangeable indices 
\begin{equation}
\label{exx}
\ex : = \{ k \in [1,N] \mid s(k) \neq + \infty \}.
\end{equation}
We will impose the following two mild conditions:

(A) The field $\KK$ contains square roots $\sqrt{\la_{kl}}$ of the scalars
$\la_{kl}$ for $1 \leq l < k \leq N$ such that the subgroup of $\kx$ generated 
by all of them contains no elements of order $2$.

(B) There exist positive integers $d_n$, $n \in \range(\eta)$,
for the function $\eta$ from Theorem \ref{chain-prime} 
such that 
$$
(\la_k^*)^{d_{\eta(l)}} = (\la_l^*)^{d_{\eta(k)}} 
$$
for all $k, l \in \ex$.

\begin{remark}
\label{A-B} With the exception of some two-cocycle twists, 
for all of the quantum nilpotent algebras $R$ coming from the 
theory of quantum groups, 
the subgroup of $\kx$ generated by $\{ \la_{kl} \mid 1 \leq l < k \leq N \}$
is torsionfree. For all such algebras, condition {\rm(A)} only requires that 
$\KK$ contains square roots of the scalars $\la_{kl}$.

All symmetric quantum nilpotent algebras that we are aware of satisfy 
$$
\la_k^* = q^{m_k} \quad {\mbox{for}} \; \; 1 \leq k \leq N
$$ 
for some non-root of unity $q \in \kx$ and positive integers 
$m_1, \ldots, m_N$. Proposition {\rm\ref{la-isi}} implies that all of them satisfy 
the condition {\rm(B)}.
\end{remark}

The set $\ex$ has cardinality $N - \rk(R)$ where $\rk(R)$ is the 
rank of the quantum nilpotent algebra $R$. By an $N \times \ex$ matrix 
we will mean a matrix of size $N \times (N - \rk(R))$ whose rows and
columns are indexed by the sets $[1,N]$ and $\ex$, respectively. The set of such
matrices with integer entries will be denoted by $M_{N \times \ex}(\Zset)$.

\begin{theorem}
\label{exchange-matr}
For every symmetric quantum nilpotent algebra $R$ 
of dimension $N$ satisfying conditions {\rm(A)} and {\rm(B)}, there exists a 
unique matrix $\wt{B} = (b_{lk}) \in M_{N \times \ex}(\Zset)$ 
whose columns satisfy the following two conditions:
\begin{enumerate}
\item[i)]
$
\Om \Big( \sum_{l=1}^N b_{lk} \ol{e_l}, \ol{e}_n \Big)
= 
\begin{cases}
\la_n^*, &\text{if} \; \; k = n
\\
1, &\text{if} \; \; k \neq n,
\end{cases}
$
for all $k \in \ex$ and $n \in [1,N]$ {\rm(}a system of linear equations written in a multiplicative form{\rm)}, and

\item[ii)]
the products $y_1^{b_{1k}} \ldots y_N^{b_{Nk}}$ 
are fixed under $\HH$ for all $k \in \ex$ {\rm(}a homogeneity condition{\rm)}.
\end{enumerate}
\end{theorem}

The second condition can be written in an explicit form using the 
fact that $y_k$ is an $\HH$-eigenvector and its eigenvalue equals the product of the $\HH$-eigenvalues of
$x_k, \ldots, x_{p^{n_k}(k)}$ where $n_k$ is the maximal nonnegative integer $n$ 
such that $p^n(k) \neq - \infty$. This property is derived from Theorem \ref{chain-prime}.

\begin{example}
\label{RqMmn-exchange-matr} It follows from Example {\rm\ref{y-quant-matr}}
that in the case of the algebras of quantum matrices $R_q[M_{m\times n}]$, 
$$
\ex = \{ (i-1)n +j \mid 1 \leq i < m, \; 1 \leq j < n \}.
$$
The bicharacter $\Om : \Zset^{mn} \times \Zset^{mn} \to \kx$ is given by 
$$
\Om( e_{(i-1)n +j}, e_{(k-1)n + l} ) = 
q^{ \de_{jl} \sign(k-i) + \de_{ik} \sign(l-j) }
$$
for all $1 \leq i,k \leq m$ and $1 \leq j,l \leq n$. Furthermore, 
$$
\la^*_s = q^2 \quad {\mbox{for}} \; \; 1 \leq s \leq mn.
$$
After an easy computation one finds that the unique solution of the 
system of equations in Theorem {\rm\ref{exchange-matr}} is given by the matrix 
$\wt{B}= (b_{(i-1)n +j, (k-1)n +j}) \in M_{mn \times \ex}(\Zset)$ with 
entries
$$
b_{(i-1)n +j, (k-1)n +l} = \begin{cases}
\pm1, &\text{if} \; \; i=k, \; l=j\pm1  \\
 &\text{or} \; \; j=l, \; k=i\pm1  \\
 &\text{or} \; \; i=k\pm1, \; j=l\pm1,  \\
0, &\text{otherwise}
\end{cases}
$$
for all $i,k \in [1,m]$ and $j,l \in [1,n]$.
\end{example}

\subsection{Cluster algebra structures}
Let us consider a symmetric quantum nilpotent algebra $R$ of dimension $N$. 
When Theorem \ref{chain-prime} is applied to the presentation of $R$
from Eq.~\eqref{tau} associated to the element $\tau \in \Xi_N$, we obtain 
a sequence of prime elements 
$$
y_{\tau, 1}, \ldots, y_{\tau, N} \in R.
$$ 
Similarly, applying Theorem \ref{exchange-matr} to the presentation 
from Eq.~\eqref{tau} we obtain the integer matrix
$$
\wt{B}_\tau \in M_{N \times \ex}(\Zset).
$$
For $\tau=\id$ we recover the original sequence $y_1, \ldots, y_N$ and 
matrix $\wt{B}$. 

\begin{theorem}
\label{main}
Every symmetric quantum nilpotent algebra $R$ of dimension $N$ satisfying 
the conditions {\rm(A)} and {\rm(B)} possesses a canonical structure 
of quantum cluster algebra for which no frozen cluster variables are inverted. 
Its initial seed has:
\begin{enumerate}
\item[i)] Cluster variables $\zeta_1 y_1, \ldots, \zeta_N y_N$ for some $\zeta_1, \ldots, \zeta_N \in \kx$,
among which the variables indexed by the set $\ex$ from Eq.~\eqref{exx} 
are exchangeable and the rest are frozen;
\item[ii)] Exchange matrix $\wt{B}$ given by Theorem {\rm\ref{exchange-matr}}.
\end{enumerate}

Furthermore, this quantum cluster algebra aways coincides with 
the corresponding upper quantum cluster algebra. 

After an appropriate rescaling,
each of the generators $x_k$ of such an algebra $R$ given by 
Eq.~\eqref{R-quant-nilp} is a cluster variable. Moreover,  
for each element $\tau$ of the subset $\Xi_N$ of the symmetric 
group $S_N$, $R$ has a seed with cluster variables obtained 
by reindexing and rescaling the sequence of prime 
elements $y_{\tau,1}, \ldots, y_{\tau, N}$. The exchange matrix 
of this seed is the matrix $\wt{B}_\tau$.
\end{theorem}

The base fields of the algebras covered by this theorem can have 
arbitrary characteristic. We refer the reader to Theorem 8.2 
in \cite{GY1} for a complete statement of the theorem, which includes 
additional results, and gives explicit formulas for the 
scalars $\zeta_1, \ldots, \zeta_N$ and the necessary reindexing and rescaling of the 
sequences $y_{\tau, 1}, \ldots, y_{\tau, N}$.

Define the following automorphism of the lattice $\Zset^N$:
$$
g = l_1 e_1 + \cdots + l_N e_N \; \mt \; \ol{g} :=l_1 \ol{e}_1 + \cdots + l_N \ol{e}_N
$$ 
for all $l_1, \ldots, l_N \in \Zset$ in terms 
of the vectors $\ol{e}_1, \ldots, \ol{e}_N$ from Eq.~\eqref{e-ol-vect}.
The construction of seeds for quantum cluster algebras in \cite{BZ}
requires assigning quantum frames to all of them. 
The quantum frame $M : \Zset^N \to \Fract(R)$ associated to the initial 
seed for the quantum cluster algebra structure in Theorem \ref{main} 
is uniquely reconstructed from the rules
\begin{align*}
&M(e_k) = \zeta_k y_k \quad {\mbox{for}} \; \; 1 \leq k \leq N \; \; \mbox{and}
\\
&M(f+g) = \Om(\ol{f}, \ol{g}) M(f) M(g) \quad {\mbox{for}} \; \; f, g \in \Zset^N.
\end{align*}
Analogous formulas describe the quantum frames associated to the 
elements $\tau$ of the set $\Xi_N$.

\begin{example}
\label{RqMmn-resc}
The cluster variables in the initial seed 
from Theorem {\rm\ref{main}}
for $R_q[M_{m \times n}]$ are 
$$
\De_{[i- \min(i,j)+1, i], [j- \min(i,j) + 1, j]}
$$
where $1 \leq i \leq m$ and $1 \leq j \leq n$. The ones 
with $i=m$ or $j=n$ are frozen. This example and Example {\rm\ref{RqMmn-exchange-matr}} 
recover the quantum cluster algebra structure of {\rm\cite{GLS}}
on $R_q[M_{m \times n}]$.
\end{example}

We finish the section by raising two questions concerning 
the line of Theorem \ref{main}:

1. If a symmetric quantum nilpotent algebra has two iterated
skew polynomial extension presentations that satisfy the assumptions 
in Definitions \ref{q-nil}, \ref{symmetric} and these two presentations are not obtained   
from each other by a permutation in $\Xi_N$, how are the corresponding 
quantum cluster algebra structures on $R$ related?

2. What is the role of the quantum seeds of a symmetric quantum nilpotent 
algebra $R$ indexed by $\Xi_N$ among the set of all quantum seeds? Is there 
a generalization of Theorem \ref{main} that constructs a larger 
family of quantum seeds using sequences of prime elements in chains
of subalgebras? 

For the first question we expect that the two quantum cluster algebra 
structures on $R$ are the same (i.e., the corresponding quantum 
seeds are mutation equivalent) if the maximal tori for the two presentations 
are the same and act in the same way on $R$. However,
proving such a fact appears to be difficult due to the generality 
of the setting. The condition on the tori is natural in light of Theorem 5.5 in 
\cite{GY0} which proves the existence of a canonical maximal torus 
for a quantum nilpotent algebra. Without imposing such a condition, the cluster 
algebra structures can be completely unrelated. For example,
every polynomial algebra
$$
R=\KK[x_1, \ldots, x_N]
$$
over an infinite field $\KK$ is a symmetric quantum nilpotent algebra 
with respect to the natural action of $(\KK^*)^N$. The quantum cluster 
algebra structure on $R$ constructed in Theorem \ref{main}
has no exchangeable indices and its frozen variables 
are $x_1, \ldots, x_N$. Each polynomial algebra has many 
different presentations associated to the elements 
of its automorphism group and the corresponding cluster 
algebra structures are not related in general.  

\section{Applications to quantum groups}

\subsection{Quantized universal enveloping algebra}
Let $\g$ be a finite dimensional complex simple Lie algebra of rank $r$  
with Cartan matrix $(c_{ij})$. For an 
arbitrary field $\KK$ and a non-root of unity $q \in \kx$, following the notation of \cite{Ja}, one defines
the quantized universal enveloping algebra $\UU_q(\g)$ 
with generators
\[
K_i^{\pm 1}, E_i, F_i, \; 1 \leq i \leq r
\]
and relations
\begin{gather*}
K_i^{-1} K_i = K_i K^{-1}_i = 1, \; \; K_i K_j = K_j K_i, 
\\
K_i E_j K^{-1}_i = q_i^{c_{ij}} E_j, \; \; 
K_i F_j K^{-1}_i = q_i^{-c_{ij}} F_j,
\\
E_i F_j - F_j E_i = \de_{i,j} \frac{K_i - K^{-1}_i}
{q_i - q^{-1}_i},
\\
\sum_{n=0}^{1-c_{ij}} (-1)^n
\begin{bmatrix} 
1-c_{ij} \\ n
\end{bmatrix}_{q_{i}}
      (E_i)^n E_j (E_i)^{1-c_{ij}-n} = 0, \; \; i \neq j,
\end{gather*}
together with the analogous relation for the generators $F_i$. 
Here $\{d_1, \ldots, d_r \}$ is the collection of 
relatively prime positive integers such that the matrix $(d_i c_{ij})$ is symmetric,
and $q_i:=q^{d_i}$. 
The algebra $\UU_q(\g)$ is a Hopf algebra 
with coproduct
\begin{gather*}
\De(K_i)   = K_i \otimes K_i,
\quad
\De(E_i) = E_i \otimes 1 + K_i \otimes E_i,
\\
\De(F_i) = F_i \otimes K_i^{-1} + 1 \otimes F_i.
\end{gather*}

The quantum Schubert cell algebras $\UU^\pm[w]$,
parametrized by the elements $w$ of the Weyl group $W$ of $\g$,
were introduced by De Concini--Kac--Procesi \cite{DKP} 
and Lusztig \cite{L}. In \cite{GLS} the term 
quantum unipotent groups was used. These algebras
are quantum analogs of the universal enveloping algebras 
$\UU(\n^\pm \cap w(\n^\mp))$ where $\n^\pm$ are the nilradicals 
of a pair of opposite Borel subalgebras of $\g$. The torus 
$\HH:= (\kx)^r$ acts on $\UU_q(\g)$ by 
\begin{equation}
\label{HUqg}
h \cdot K_i^{\pm 1} = K_i^{\pm 1}, \; \;
h \cdot E_i = \xi_i E_i, \; \;
h \cdot F_i = \xi_i^{-1} F_i,
\end{equation} 
for all $h = (\xi_1, \ldots, \xi_r) \in \HH$ and $1 \leq i \leq r$.
The subalgebras $\UU^\pm[w]$ are preserved by this action.

Denote by $\al_1, \ldots, \al_r$ the set of simple roots 
of $\g$ and by $s_1, \ldots, s_r \in W$ the 
corresponding set of simple reflections.
All algebras $\UU^\pm[w]$ are symmetric quantum nilpotent 
algebras for all base fields $\KK$ and non-roots of unity $q \in \kx$.
In fact, to each reduced expression
$$
w = s_{i_1} \ldots s_{i_N},
$$
one associates a presentation of $\UU^\pm[w]$ that 
satisfies Definitions \ref{q-nil} and \ref{symmetric} 
as follows. Consider the Weyl group elements
$$
w_{\leq k} := s_{i_1} \ldots s_{i_k}, \; 1 \leq k \leq N, \quad \mbox{and} 
\quad w_{\leq 0 }:=1.
$$
In terms of these elements the roots of the Lie algebra $\n^+ \cap w(\n^-)$
are 
$$
\be_1 = \al_{i_1}, \; \be_2 = w_{\leq 1} (\al_{i_2}), \; \ldots, \;
\beta_N := w_{\leq N-1} (\al_{i_N}).
$$
As in \cite{Ja,L}, one associates with those roots
the Lusztig root vectors $E_{\be_1}, F_{\be_1}, \ldots, E_{\be_N}, F_{\be_N} \in \UU_q(\g)$. 
The quantum Schubert cell algebra $\UU^-[w]$, defined as the subalgebra of $\UU_q(\g)$ generated by $F_{\be_1}, \ldots, F_{\be_N}$, has an iterated skew polynomial 
extension presentation of the form
\begin{equation}
\label{UwOre}
\UU^-[w] = \KK [F_{\be_1}] [F_{\be_2}; \sig_2, \de_2] 
\ldots [F_{\be_N}; \sig_N, \de_N]  
\end{equation}
for which conditions (a)--(c) of Definition \ref{q-nil} are satisfied with respect to the action 
Eq.~\eqref{HUqg}. Moreover, this presentation satisfies the condition for a symmetric 
quantum nilpotent algebra in Definition \ref{symmetric} because 
of the Levendorskii--Soibelman straightening law \cite{Ja} in $\UU_q(\g)$.
The opposite algebra $\UU^+[w]$, generated by $E_{\be_1}, \ldots, E_{\be_N}$,
has analogous properties and is actually isomorphic \cite{Ja} to $\UU^-[w]$.

The algebra $R_q[M_{m\times n}]$ is isomorphic to one of the algebras
$\UU^-[w]$ for $\g = {\mathfrak{sl}}_{m+n}$ and a certain choice of $w \in S_{m+n}$.

\subsection{Quantized function algebras}
\  The irreducible finite dimensional modules of $\UU_q(\g)$ on which the elements $K_i$ act 
diagonally via powers of  the scalars $q_i$ are parametrized by the set $P_+$ of dominant integral weights of
$\g$. The module corresponding to such a weight $\la$ will be denoted by $V(\la)$. 

Let $G$ be the connected, simply connected algebraic group with Lie algebra $\g$.
The Hopf subalgebra of $\UU_q(\g)^*$ spanned by the matrix coefficients $c^\la_{f, y}$
of all modules $V(\la)$ (where $f \in V(\la)^*$ and $y \in V(\la)$)
is denoted by $R_q[G]$ and called the quantum group corresponding to $G$.
The weight spaces $V(\la)_{w \la}$ are one dimensional for all Weyl group
elements $w$. Considering the fundamental representations $V(\vpi_1), \ldots, V(\vpi_r)$, 
and a normalized covector and vector in each of those weight spaces, one defines, following \cite{BZ},
the quantum minors 
$$
\De^i_{w, v} = \De_{w \vpi_i, v \vpi_i} \in R_q[G], \quad 1 \leq i \leq r, \; w, v \in W.
$$   
The subalgebras of $R_q[G]$ spanned by the elements of the form $c_{f,y}^\la$ 
where $y$ is a highest or a lowest weight vector of $V(\la)$ and $f \in V(\la)^*$ 
are denoted by $R^\pm$. They are quantum analogs of the base affine space
of $G$. With the help of the Demazure 
modules $V_w^+(\la) = \UU^+_q(\g) V(\la)_{w \la}$ 
(where $\UU^+_q(\g)$ is the unital subalgebra of $\UU_q(\g)$ generated by 
$E_1, \ldots, E_r$) one defines the ideal 
$$
I^+_w = \Span \{ c^\la_{f, y} \mid \la \in P_+, \; f|_{V^+_w(\la)}= 0, \; y \in V(\la)_\la \} 
$$
of $R^+$. Analogously one defines an ideal $I^-_w$ of $R^-$. 
For all pairs of Weyl group elements $(w,v)$,
the quantized coordinate ring of the double Bruhat cell \cite{FZ1}
$$
G^{w,v}: = B^+ w B^+ \cap B^- v B^- ,
$$
where $B^\pm$ are opposite Borel subgroups of $G$, is defined by 
$$
R_q[G^{w,v}] := (I^+_w R^- + R^+ I^-_v) [ (\De^i_{w, 1})^{-1}, 
(\De^i_{v w_0, w_0})^{-1}]
$$
where the localization is taken over all $1 \leq i \leq r$ and 
$w_0$ denotes the longest element of the Weyl group $W$. 

To connect cluster algebra structures on the two kinds of
algebras (quantum Schubert cells and quantum double Bruhat cells), 
we use certain subalgebras $S^+_w$ 
of $(R^+/I_w^+)[(\De^i_{w, 1})^{-1}, 1 \leq i \leq r]$, which were
defined by Joseph \cite{J}.
They are the subalgebras generated by the elements 
$$
(\De^i_{w, 1})^{-1} (c^{\vpi_i}_{f, y} + I_w^+),
$$
for $1 \leq i \leq r$, $f \in V(\vpi_i)^*$, 
and $y$ a highest weight vector of $V(\vpi_i)$.
These algebras played a major role in the study of the 
spectra of quantum groups \cite{J,Y-mem} 
and the quantum Schubert cell algebras \cite{Y-pmls}. In \cite{Y-mem},
Yakimov constructed an algebra antiisomorphism
$$
\vp_w : S^+_w \to \UU^-[w].
$$
An earlier variant of it for $\UU_q(\g)$ equipped with 
a different coproduct appeared in \cite{Y-pmls}.

\subsection{Cluster structures on quantum Schubert cell algebras}
Denote by $\lcor.,. \rcor$ the Weyl group invariant bilinear form 
on the vector space $\Rset \al_1 \oplus \cdots \oplus \Rset \al_r$
normalized by $\lcor \al_i, \al_i \rcor = 2$ for short roots $\al_i$.
Let $\|\ga \|^2:= \lcor \ga, \ga \rcor$. 

Fix a Weyl group element $w$ 
and consider the quantum Schubert cell algebra $\UU^-[w]$.
A reduced expression 
$w =s_{i_1} \ldots s_{i_N}$ gives rise to the presentation in 
Eq.~\eqref{UwOre} of $\UU^-[w]$ as a symmetric 
quantum nilpotent algebra. The result of the application 
of Theorem \ref{chain-prime} to it is as follows. The function $\eta$ 
can be chosen as 
$$
\eta(k) = i_k \quad \mbox{for all} \; \; 1 \leq k \leq N. 
$$
The predecessor function $p$ is the function $k \mt k^-$ 
which plays a key role in the works of Fomin and Zelevinsky \cite{FZ2,FZ1}, 
$$
k^-
:=
\begin{cases}
\max \{ l <k \mid i_l = i_k \}, 
&\text{if such $l$ exists}, 
\\
- \infty, \; & \text{otherwise}. 
\end{cases}
$$
The successor function $s$ is the function $k \mt k^+$ in \cite{FZ2,FZ1}.
The sequence of prime elements $y_1, \ldots, y_N$ consists of scalar multiples of the elements
$$
\vp_w \left( \De_{w_{\leq p^{n_k}(k)-1}, 1}^{i_k} (\De_{w_{\leq k}, 1}^{i_k})^{-1} \right), 
\quad 1 \leq k \leq N
$$
where $n_k$ denotes the maximal nonnegative integer $n$ such that $p^n(k) \neq - \infty$. 
The presentation in Eq.~\eqref{UwOre} of $\UU^-[w]$ as a symmetric quantum nilpotent algebra 
satisfies the conditions 
(A) and (B) if $\sqrt{q} \in \KK$, in which case Theorem \ref{main} produces a canonical 
cluster algebra structure on $\UU^-[w]$. Among all the clusters in Theorem \ref{main}, 
indexed by the elements of the subset $\Xi_N$ of the symmetric group $S_N$, the 
one corresponding to the longest element of $S_N$ is closest to the combinatorial 
setting of \cite{BFZ,BZ}. It goes with the reverse presentation of $\UU^-[w]$,
$$
\UU^-[w]= \KK [F_{\be_N}] [F_{\be_{N-1}}; \sig^*_{N-1}, \de^*_{N-1}] \ldots [F_{\be_1}; \sig^*_1, \de^*_1].
$$
The transition from the original presentation in Eq.~\eqref{UwOre} 
to the above one amounts to interchanging the roles of the predecessor 
and successor functions. As a result, the set of exchangeable 
indices for 
the latter presentation is 
\begin{equation}
\label{ex-w}
{\ex}_w := \{ k \in [1,N] \mid k^- \neq - \infty \}.
\end{equation} 

\begin{theorem}
\label{Uw-thm} Consider an arbitrary finite dimensional complex simple Lie algebra $\g$, a Weyl group element $w \in W$,
a reduced expression of $w$, an arbitrary base field $\KK$ and a non-root of 
unity $q \in \kx$ such that $\sqrt{q} \in \KK$. The quantum Schubert cell algebra 
$\UU^-[w]$ possesses a canonical quantum cluster algebra structure for which 
no frozen cluster variables are inverted and the
set of exchangeable indices is $\ex_w$. Its
initial seed consists of the cluster variables 
$$
{\sqrt{q}}^{\, \| (w - w_{\leq k-1})\vpi_{i_k} \|^2/2 }
\vp_w \left( \De_{w_{\leq k-1}, 1}^{i_k} (\De_{w, 1}^{i_k})^{-1} \right),
$$
$1 \leq k \leq N$. 
The exchange matrix $\wt{B}$ of this seed has entries given by 
$$
b_{kl} = 
\begin{cases}
1, &\text{if} \; \; k = p(l)
\\
-1, &\text{if} \; \; k = s(l)
\\
c_{i_k i_l}, & \text{if} \; \; 
p(k) < p(l) < k < l
\\ 
- c_{i_k i_l}, & \text{if} \; \; 
p(l) < p(k) < l < k
\\
0, & \text{otherwise}
\end{cases}
$$
for all $1 \leq k \leq N$ and $l \in {\ex}_w$. Furthermore, this 
quantum cluster algebra equals the corresponding 
upper quantum cluster algebra. For all $k \in [1,N]$, 
$m \in \Zset_{\geq 0}$ such that $s^m(k) \in [1,N]$, the 
elements  
\begin{multline*}
{\sqrt{q}}^{\, \| (w_{\leq s^m(k)} - w_{\leq k-1})\vpi_{i_k} \|^2/2 } \times 
\\
\times \vp_{w_{\leq s^m(k)}} \left( \De_{w_{\leq k-1}, 1}^{i_k} (\De_{w_{\leq s^m(k)}, 1}^{i_k})^{-1} \right)
\end{multline*}
are cluster variables of $\UU^-[w]$.
\end{theorem}

For symmetric Kac-Moody algebras $\g$ the theorem is due to Gei\ss, Leclerc, 
and Schr\"oer \cite{GLS}. Our proof also works for all Kac-Moody algebras $\g$, 
but here we restrict to the finite dimensional case for simplicity of the exposition.
Theorem \ref{Uw-thm} is proved in Section 10 of \cite{GY1}.

Examples \ref{RqMmn-exchange-matr} and \ref{RqMmn-resc} can be recovered as special cases of Theorem \ref{Uw-thm} 
for $\g = \sl_{m+n}$ and a particular choice of the Weyl group element $w \in S_{m+n}$. 
In this case the torus action can be used to kill the power of $\sqrt{q}$.

\begin{remark}
\label{cluster} It follows from the definition of the antiisomorphism
$\vp_w : S_w^+ \to \UU^-[w]$ in {\rm\cite{Y-mem,Y-pmls}}  that the element 
\begin{equation}
\label{element}
\vp_w \left( \De_{w_{\leq k-1}, 1}^{i_k} (\De_{w, 1}^{i_k})^{-1} \right)
\end{equation}
is obtained (up to a minor term) by evaluating 
$$
\left( \De_{w_{\leq k-1}, w}^{i_k} \otimes \id \right) ( {\mathcal{R}}^w )
$$
where ${\mathcal{R}}^w$, called the $R$-matrix for the Weyl group element $w$, 
equals the infinite sum $\sum_j u_j^+ \otimes u^-_j$ for dual bases $\{u_j^+\}$ and $\{u^-_j\}$ 
of $\UU^+[w]$ and $\UU^-[w]$. Because of this, the element 
in Eq.~\eqref{element} can be identified with 
$\De_{w_{\leq k-1}, w}^{i_k}$ 
and thus can be thought of as a quantum minor. 
Such a construction of cluster variables of $\UU^-[w]$ 
via quantum minors is due to {\rm\cite{GLS}}, which used linear maps that are 
not algebra {\rm(}anti\/{\rm)}isomorphisms.

More generally, 
$$
\UU^-[w_{\leq j}] \subseteq \UU^-[w] \quad {\mbox{for}} \; \; 1 \leq j \leq N
$$
and the cluster variables in Theorem {\rm\ref{Uw-thm}},
$$
\vp_{w_{\leq s^m(k)}} \left( \De_{w_{\leq k-1}, 1}^{i_k} (\De_{w_{\leq s^m(k)}, 1}^{i_k})^{-1} \right)
\in \UU^-[w_{\leq s^m(k)}],
$$
can be identified with the quantum minors $\De_{w_{\leq k-1}, w_{\leq s^m(k)}}^{i_k}$.
\end{remark}

\subsection{The Berenstein--Zelevinsky conjecture}
\  Consider a pair of Weyl group elements $(w,v)$ with 
reduced expressions 
$$
w=s_{i_1} \ldots s_{i_N} \quad \text{and} \quad 
v = s_{i'_1} \ldots s_{i'_M}.
$$
Let $\eta : [1,r + M + N] \to [1,r]$ be the function given by
$$
\eta(k) = \begin{cases}
k, &\text{for} \; \; 1 \leq k \leq r,
\\
i'_{k-r}, &\text{for}\; \; r+1 \leq k \leq r+M, 
\\
i_{k-r-M}, &\text{for} \; \; 
r+M+1 \leq k \leq r+M + N.
\end{cases}
$$
The following set will be used as the set of exchangeable indices 
for a quantum cluster algebra structure on $R_q[G^{w,v}]$:
$$
\ex_{w,v}:= [1,k] \sqcup \{ k \in [r+1, r+M+N] \mid 
s(k) \ne + \infty \}
$$
where $s$ is the successor function for the level sets of $\eta$.
Set $\ep(k):=1$ for $k \leq r+M$ and $\ep(k) := - 1$ for $k >r+M$.
Following \cite{BZ}, define the $(r+M+N) \times \ex$ matrix $\wt{B}_{w,v}$ 
with entries
$$
b_{kl} := 
\begin{cases}
-\ep(l), &\text{if} \; \; k = p(l), 
\\
- \ep(l) c_{\eta(k), \eta(l)}, 
&\text{if $k < l < s(k) < s(l), \, \ep(l) = \ep(s(k))$}
\\
&\text{or $k < l \leq r+M < s(l) < s(k)$,}
\\
\ep(k) c_{\eta(k), \eta(l)}, 
&\text{if $l < k < s(l) < s(k), \, \ep(k) = \ep(s(l))$}  
\\
&\text{or $l < k \leq r+M < s(k) < s(l)$,}
\\
\ep(k), &\text{if $k = s(l)$,}
\\
0, & \text{otherwise}.
\end{cases}
$$

The following theorem proves the Berenstein--Zelevinsky conjecture, \cite{BZ}.

\begin{theorem}
\label{BZconj}
Let $G$ be an arbitrary complex simple Lie group and 
$(w,v)$ a pair of elements of the corresponding Weyl group. 
For any base field $\KK$ and a non-root of unity $q \in \kx$ 
such that $\sqrt{q} \in \kx$, the quantum double Bruhat cell 
algebra $R_q[G^{w,v}]$ possesses a canonical structure 
of quantum cluster algebra for which all frozen cluster variables are inverted and 
the set of exchangeable indices is $\ex_{w,v}$.  
The initial seed has exchange matrix $\wt{B}_{w,v}$ defined above
and cluster variables $y_1, \ldots, y_{r+M+N} \in R_q[G^{w,v}]$ 
given by 
$$
y_k = 
\begin{cases}
\De^k_{1,v^{-1}}, &\text{for} \; \; 1 \leq k \leq r,
\\
\xi_k
\De^{i'_{k-r}}_{1, v^{-1} v_{\leq k-r}}, 
&\text{for}\; \; r+1 \leq k \leq r+M,
\\
\xi_k
\De^{i_{k-r-M}}_{w_{\leq k -r -M}, 1}, 
&\text{for} \; \; 
r+M+1 \leq k \leq r+M + N
\end{cases}
$$
for some scalars $\xi_k \in \kx$.

Furthermore, this quantum cluster algebra coincides with the corresponding upper 
quantum cluster algebra.
\end{theorem}

We briefly sketch the relationship of the quantum double Bruhat cell algebras $R_q[G^{w,v}]$ 
to quantum nilpotent algebras and the proof of the theorem. We first show, using results of Joseph \cite{J}, 
that $R_q[G^{w,v}]$ is a localization of 
$$
(S^+_w \bowtie S^-_v) \,\#\, \KK[(\De^1_{1, v^{-1}})^{\pm 1} , \ldots, (\De^r_{1, v^{-1}})^{\pm 1}]
$$
where $S^-_v$ is the Joseph subalgebra of $R^-$ defined
in a similar way \cite{J} to the subalgebra $S^+_w$ of $R^+$. The ``bicrossed'' and 
smash products are defined \cite{J,Y-mem} from the Drinfeld $R$-matrix commutation relations  
of $R_q[G]$. Using the antiisomorphism $\vp_w$ and its negative counterpart 
(which turns out to be an isomorphism \cite{Y-mem}), 
one converts 
$$
S^+_w \bowtie S^-_v \cong {\UU^-[w]}^{\mathrm{op}} \bowtie {\UU^+[v]},
$$
where $R^{\mathrm{op}}$ stands for the algebra with opposite product. We then establish 
that the right hand algebra above is a symmetric quantum nilpotent algebra 
satisfying the conditions (A) and (B). The proof is completed by applying 
Theorem \ref{main}, and showing that the localization that we started with is, in fact, 
a localization by all frozen cluster variables.

\section{Poisson nilpotent algebras and cluster algebras}
In this subsection, we will assume that the base field $\KK$ has characteristic $0$.
A prime element $p$ of a Poisson algebra $R$ with Poisson bracket $\{., .\}$ will be 
called Poisson prime if 
$$
\{R, p \} = Rp.
$$ 
In other words, this requires that the principal ideal $Rp$ be a Poisson ideal as well as a prime ideal.

For a commutative algebra $R$ equipped with a rational action of a 
torus $\HH$ by algebra automorphisms, we will denote by 
$\partial_h$ the derivation of $R$ corresponding to an element $h$ of 
the Lie algebra of $\HH$. 

\begin{definition}
\label{Pnilpotent}
A nilpotent semi-quadratic Poisson algebra is a polynomial algebra $\KK[x_1, \ldots, x_N]$ 
with a Poisson structure $\{.,.\}$ and a rational action of a torus 
$\HH=(\kx)^r$ by Poisson algebra automorphisms for which $x_1, \ldots, x_N$ are 
$\HH$-eigenvectors 
and there exist
elements $h_1, \ldots, h_N$ in the Lie algebra of $\HH$ such that
the following two conditions are satisfied for $1 \leq k \leq N$:

{\rm(a)} For all $b \in R_{k-1}:=\KK[x_1, \ldots, x_{k-1}]$
$$
\{x_k, b \} = \partial_{h_k}(b) x_k + \de_k(b)
$$ 
for some $\de_k(b) \in R_{k-1}$ and the map $\de_k \colon R_{k-1} \to R_{k-1}$
is locally nilpotent.

{\rm(b)} The $h_k$-eigenvalue of $x_k$ is non-zero.
\\
Such an algebra will be called symmetric if the 
above condition is satisfied for the reverse order of generators $x_N, \ldots, x_1$
{\rm(}with different choices of elements $h_\bullet${\rm)}.  
\end{definition}

The adjective semi-quadratic refers to the leading term in the Poisson bracket 
$\{x_k, x_l\}$ which is forced to have the form $\la_{kl} x_k x_l$ for some 
$\la_{kl} \in \kx$.

\begin{theorem}
\label{Poisson} Every symmetric nilpotent semi-quadratic Poisson algebra as above satisfying 
the Poisson analog of condition {\rm(B)} has a canonical structure of cluster algebra 
for which no frozen variables are inverted and 
the compatible Poisson bracket in the sense of Gekhtman--Shapiro--Vainshtein
{\rm\cite{GSV}} is $\{.,.\}$. Its initial cluster consists, up to scalar multiples, of those Poisson prime elements of the 
chain of Poisson subalgebras
$$
\KK[x_1] \subset \KK [x_1, x_2] \subset \ldots \subset \KK[x_1, \ldots, x_N]
$$   
which are $\HH$-eigenvectors. Each generator $x_k$, $1 \leq k \leq N$, 
is a cluster variable of this cluster algebra.

Furthermore, this cluster algebra coincides with the corresponding 
upper cluster algebra.
\end{theorem}

Applying this theorem in a similar fashion to Theorem \ref{BZconj} we
obtain the following result.

\begin{theorem} 
\label{doubleBruhat}
Let $G$ be an arbitrary complex simple Lie group. For all pairs of elements $(w,v)$
of the Weyl group of $G$, the Berenstein--Fomin--Zelevinsky upper cluster algebra 
{\rm\cite{BFZ}} on the coordinate ring of the double Bruhat cell $G^{w,v}$ coincides with 
the corresponding cluster algebra. In other words, the coordinate rings of all 
double Bruhat cells $\Cset[G^{w,v}]$ are cluster algebras with initial seeds 
constructed in {\rm\cite{BFZ}}.
\end{theorem}

\begin{acknowledgments} 
We are grateful to S. Fomin, A. Berenstein, Ph. Di Francesco, R. Kedem, B. Keller, 
A. Knutson, 
B. Leclerc, N. Reshetikhin, D. Rupel and 
A. Zelevinsky for helpful discussions and comments. 
We are also thankful to the referee, whose suggestions 
helped improve the exposition.
Moreover, we would like to thank
the Mathematical Sciences Research Institute for its hospitality during the programs in ``Cluster Algebras'' 
and ``Noncommutative Algebraic Geometry and Representation 
Theory'' when parts of this project were completed.
This work was partially supported by the National Science 
Foundation grants DMS-0800948 (to K.R.G.), and DMS-1001632 and DMS-1303038 (to M.T.Y.).  
\end{acknowledgments}

\end{article}
\end{document}